\documentclass{amsart}
\usepackage{amsmath,amscd}
\usepackage{enumerate}
\usepackage[dvips]{graphicx}
\usepackage{rotating}
\usepackage{cite}
\usepackage{amsmath,amsfonts}
\usepackage{latexsym,amssymb,amscd}
\newtheorem{thm}{Theorem}[section]
\newtheorem{cor}[thm]{Corollary}
\newtheorem{lem}[thm]{Lemma}
\newtheorem{defn}[thm]{Definition}
\newtheorem{rem}[thm]{Remark}
\newtheorem{notn}[thm]{Notation}
\newtheorem{ex}[thm]{Example}
\newtheorem{alg}[thm]{Algorithm}

\newtheorem{prob}[thm]{Problem}
\def\sqr#1#2{{\vcenter{\hrule height.#2pt
        \hbox{\vrule width.#2pt height#1pt \kern#1pt
                \vrule width.#2pt}
        \hrule height.#2pt}}}
\def\square{\mathchoice\sqr64\sqr64\sqr{4}3\sqr{3}3}
\def\QED{\hfill$\square$\bigskip}
\begin{document}
\title[Some Bounds in the B. and M. Shapiro Conjecture for Flag Varieties]{Some Lower Bounds in the B. and M. Shapiro Conjecture for Flag Varieties}

\author{Monique Azar}
\address{Department of Mathematics,
American University of Beirut, Beirut, Lebanon}
\email{monique.azar@aub.edu.lb}
\author{Andrei Gabrielov}
\address{Department of Mathematics,
Purdue University, West Lafayette, IN 47907, USA}
\email{agabriel@math.purdue.edu}
\thanks{The second  author was supported by the NSF grant DMS-0801050.}
\begin{abstract}
The B.\ and M.\ Shapiro conjecture stated that all solutions
of the Schubert Calculus problems associated with 
real points on the rational normal curve should be real.
For Grassmannians, it was proved by Mukhin, Tarasov and Varchenko.
For flag varieties, Sottile found a counterexample
and suggested that all solutions should be real
under certain monotonicity conditions.
In this paper, we compute lower bounds on the number of real
solutions for some special cases of the  B.\ and M.\ Shapiro 
conjecture for flag varieties, when Sottile's monotonicity 
conditions are not satisfied.
\end{abstract}
\maketitle

\section{Introduction}
Schubert Calculus is a recipe for counting
geometric objects subject to certain incidence relations
~\cite{SC}. For example,

\begin{prob}
	\label{Clines} Given $2d-2$ lines in general position in 
	$\mathbb{CP}^d$. How many codimension 2 subspaces of 
	$\mathbb{CP}^d$ meet all $2d-2$ lines?
\end{prob}

The answer is
$u_d=\frac{1}{d}{{2d-2}\choose{d-1}}$, the $d$-th Catalan number. 
Schubert Calculus was based on the semi-empirical principle of conservation of 
number. Its rigorous foundation,
subject of Hilbert's 15-th problem, was established 
through the development of intersection theory. 
>From the point of view of intersection theory, enumerative 
problems such as Problem \ref{Clines} are solved by counting intersection 
multiplicities of Schubert varieties in the Grassmannian (see, e.g.,~\cite{B}).

Over the reals, the principle of conservation of number fails,
and solving enumerative problems becomes considerably more complicated.
In~\cite{F}, W. Fulton says
``The question of how many solutions of real equations can be 
real is still very much open, particularly for enumerative 
problems.'' 

For example, suppose that all $2d-2$ lines in Problem \ref{Clines} are real.
Of all the codimension 2 subspaces of $\mathbb{CP}^d$ 
that meet all these lines, how many are real?
Sottile~\cite{So1} proved that all of them can be real for some
choice of the $2d-2$ lines.
Boris and Michael Shapiro conjectured that all these subspaces are 
real whenever the given lines are tangent to the rational normal curve
$\gamma(t)=[1:t:t^2:\dots :t^d]$ at distinct real points. Eremenko and 
Gabrielov~\cite{EG1} proved the following equivalent theorem.

\begin{thm}\label{eg}
	If the critical points of a rational function $f$ are all real then 
	$f$ is equivalent to a real rational function, i.e., there exists a 
	fractional linear transformation $l$ such that $l\circ f$ is a real 
	rational function.
\end{thm}

A more general version of the B.\ and M.\ Shapiro conjecture, proved in
~\cite{MTV1}\footnote{This result dates back to 2005.}, claims a similar result for higher dimensional subspaces.

\begin{thm} 
	Let $P_1, P_2, \ldots, P_{k(d-k+1)}$ be $(k-1)$-dimensional planes in 
	$\mathbb{CP}^d$ osculating the rational normal curve at distinct real 
	points. Any $(d-k)$-dimensional plane that meets all 
	$P_i$, $1\leq i\leq k(d-k+1)$, must be real.
\end{thm}

B. and M. Shapiro suggested an extension of their conjecture
to flag varieties, replacing osculating planes 
by osculating flags. A special case of that conjecture
would imply that all solutions to the following
problem are real.

\begin{prob}\label{onesecant} 
	Let $y_1,y_2,\ldots , y_{2d-3},r,s$ be distinct real points. For $1\leq i\leq 2d-3$, 
	let $T_i$ be the line tangent to $\gamma$ at $\gamma(y_i)$ and let $T_{2d-2}$ be 
	the line through the points $\gamma(r)$ and $\gamma(s)$. Among all 
	the codimension 2 subspaces of $\mathbb{CP}^d$ that meet all the lines 
	$T_1,T_2,\ldots , T_{2d-2}$, how many are real?
\end{prob}

Sottile~\cite{So1} found examples when the B. and M. Shapiro conjecture
for flag varieties fails, in particular, Problem \ref{onesecant}
does have non-real solutions.
Computer experiments suggested that the conjecture might hold whenever a 
certain monotonicity condition is met~\cite{RSSS}. 
For the case of Problem \ref{onesecant}, Sottile's 
monotonicity condition simply means that the interval $I$ with endpoints 
$r$ and $s$ contains either all or none of the points 
$y_1,y_2,\ldots , y_{2d-3}$. A proof of this special case was given in~\cite{EGSV}.
For a complete survey of the B. and M. Shapiro conjecture and its modifications, see~\cite{So2} or~\cite{So3}.

In this paper we study Problem \ref{onesecant} when the monotonicity 
condition does not hold, i.e., when the interval $I$ contains $k$
of the points $y_1, y_2,\ldots ,  y_{2d-3}$ 
for some $k$ with $1\leq k<2d-3$. We give an algorithm 
to compute a (non-strict)
lower bound for the number of real subspaces for any pair of 
integers $d$ and $k$ with $1\leq k< 2d-3$. We do this by giving a lower bound for the number of solutions to the following equivalent problem.

\begin{prob}\label{rrf} 
	Let $y_1, y_2,\ldots , y_{2d-3},r,s$ be distinct real points and suppose that the interval $I$ having endpoints $r$ and   $s$ contains $k$ of the points $y_1, y_2,\ldots, y_{2d-3}$ for some $k\in\{1,2,\ldots,2d-4\}$. How many equivalence classes of real rational functions $f$ of degree $d$ having critical points at 
	$y_1,y_2,\ldots, y_{2d-3}$ and satisfying $f(r)=f(s)$ are there? 
\end{prob}

\begin{rem}\label{improvelb} The answer to Problem \ref{rrf} does not change if we replace $k$ by $2d-3-k$.\end{rem}

The lower bound for the number of real solutions of Problem \ref{rrf}
is obtained by studying a one-parametric family 
of real rational functions with critical points $y_1,\ldots, y_{2d-3},\theta$,
and the dependence on $\theta$ of increments on $[r,s]$ 
of the functions in that family.

\section{The Wronski map and one-parametric families of rational functions.}

In this section, we represent the points of the
Grassmannian $G$ of 2-dimensional planes
in the space of polynomials of degree $d$ by pencils $(p,q)$
of linearly independent polynomials,
define the Wronski map  $\mathcal{W}:G\longrightarrow\mathbb{P}^{2d-2}$
and prove that $\mathcal{W}$ is not ramified over the space $Q$ of
polynomials with all real roots of multiplicity at most 2.
We study the properties of a one-parametric family $f_t=p_t/q_t$
of rational functions obtained by lifting to $G$ a path
$t\mapsto w_{\theta(t)}$ in $Q$.
Here $\theta(t)$ is a path in $\mathbb{\overline{R}}$ and
$w_{\theta}$ is a polynomial with roots $y_1,\dots,y_{2d-3},\theta$.
Our main tool is a net of a real rational function
with real critical points (see \cite{EG1}).

Let $G$ be the Grassmannian of 2-dimensional planes in 
the space of complex polynomials of degree at most $d$. 
An element of $G$ can be defined by a pencil $L=L(p,q)=\{\alpha p+\beta q\}$ 
where $p$ and $q$ are two linearly independent polynomials of degree at most $d$ 
and $[\alpha:\beta]\in\mathbb{CP}^1$. 
If $L\in G_{\mathbb{R}}$, the Grassmannian of real two-dimensional planes, 
we always assume that $p$ and $q$ are real and call $L$ a real pencil. 
The Wronskian of $(p,q)$, given by $W(p,q)=p'q-q'p$, is a nonzero polynomial of 
degree at most $2d-2$. 
If the degree of $W(p,q)$ is $2d-2-k$ for some $k>0$, we say that $W(p,q)$ 
has a root at infinity of multiplicity $k$. 
If $(p_1,q_1)$ is another basis for $L$, then $W(p_1,q_1)$ differs from
$W(p,q)$ by a nonzero multiplicative constant. 
For a nonzero polynomial $r=\sum a_iz^i$ of degree at most $2d-2$, 
let $h(r)=[a_0:a_1:\ldots:a_{2d-2}]$. 

\begin{defn} The Wronski map $\mathcal{W}:G\longrightarrow\mathbb{P}^{2d-2}$ is defined by 
$\mathcal{W}(L)=h(W(p,q))$, where $(p,q)$ is a basis for $L$.
\end{defn}

The map $\mathcal{W}$ is well-defined and finite. 
Its degree is $u(d)$, the $d$-th Catalan number. 
This was originally computed by Schubert in 1886. 
A proof can be found in~\cite{Go} or~\cite{HP}. 

\begin{defn}Let $Q=h(Q_0)$, where $Q_0$ is the set of all polynomials $r$ satisfying

{\rm (i)} $2d-4\le\deg(r)\leq 2d-2$,

{\rm (ii)} all roots of $r$ belong to $\mathbb{\overline{R}}$, and

{\rm (iii)} all roots of $r$ have multiplicity at most 2.
 
\end{defn}

In Theorem \ref{unram}, we shall show that $\mathcal{W}$ is unramified over $Q$ and hence, 
$\mathcal{W}|_{\mathcal{W}^{-1}(Q)}$ is a covering map. 
Thus, for any path $\rho$ in $Q$ with initial point $\rho_0$, given $L_0\in G$ with 
$\mathcal{W}(L_0)=\rho_0$, $\rho$ can be lifted in a unique way to a path in $G$ 
starting at $L_0$.

\begin{defn}Let $D=\{y_1,y_2,\ldots,y_{2d-3}\}$ be a set of points in $\mathbb{R}$. 
For $y\in\mathbb{R}$, let $\omega_y=(x-y_1)(x-y_2)\ldots (x-y_{2d-3})(x-y)$, and let 
$\omega_\infty=(x-y_1)(x-y_2)\ldots (x-y_{2d-3})$.
\end{defn}
 
The map $\mathbb{\overline{R}}\longrightarrow \mathbb{P}^{2d-2}$ given by 
$y \longrightarrow h(w_y)$ is continuous. 
In fact, it extends to an algebraic map 
$\mathbb{\overline{C}}\longrightarrow \mathbb{CP}^{2d-2}$.
Let $\mathcal{I}$ be an interval in $\mathbb{R}$ and let $t_0\in \mathcal{I}$. 
A path $\theta:\mathcal{I}\longrightarrow \mathbb{\overline{R}}$ induces a path 
$\rho:\mathcal{I}\longrightarrow Q$ given by $\rho(t)=h(\omega_{\theta(t)})$. 
Given $L_0\in G$ with $\mathcal{W}(L_0)=h(\omega_{\theta(t_0)})$, $\rho$ can be lifted 
to a unique path $\chi$ in $G$ satisfying $\chi(t_0)=L_0$. 
It was shown in ~\cite{EG1} that $\mathcal{W}^{-1}(\mathbb{RP}^{2d-2})\subset G_\mathbb{R}$. 
This implies that $\chi(\mathcal{I})$ is contained in $G_{\mathbb{R}}$. 

\begin{lem}\label{lifts} Let $\mathcal{R}_d$ be the space of real rational functions of 
degree at most $d$ all of whose critical points are real. 
If $\theta$ is analytic, all the lifts defined above are analytic, and polynomials 
$p_t$ and $q_t$ with coefficients analytic in $t$ can be selected to form the bases 
of $\chi(t)$, $t\in\mathcal{I}$. 
Let $f_t$ be the real rational function $p_t/q_t$ viewed as a holomorphic function 
from $\mathbb{\overline{C}}$ to $\mathbb{\overline{C}}$. 
The path $\eta:\mathcal{I}\longrightarrow \mathcal{R}_d$ given by 
$t\longrightarrow f_t$ has the following properties.

{\rm(i)} The Wronskian  $W(p_t,q_t)$ is equal to $\rho(t)$. 

{\rm(ii)} If $x\in \mathbb{\overline{C}}\setminus D$ then the map 
$t\stackrel{\phi_x}{\longrightarrow} f_t(x)$ is continuous. 

{\rm(iii)} For each $x\in D$, the map  $t\stackrel{\phi_x}{\longrightarrow} f_t(x)$ 
is continuous at all values $t$ for which $p_t$ and $q_t$ do not have a common root at $x$. 

\end{lem}

\textbf{Proof. }  
The first property follows from the construction of $f_t$. 
To prove (ii) and (iii), first observe that if $x\in\mathbb{\overline{C}}\setminus D$ then, 
for any $t\in\mathcal{I}$, $p_t(x)$ and $q_t(x)$ cannot both be zero. 
Let $(t_0,x)\in\mathcal{I}\times\mathbb{\overline{C}}$ be such that at least one of 
$p_{t_0}(x)$ and $q_{t_0}(x)$ is not equal to $0$. 
This implies that there exists a neighborhood $N$ of $t_0$ such that either 
$f_{t}(x)$ or $1/f_{t}(x)$ is a continuous complex-valued function on $N$ and hence, 
$f_t(x):N\longrightarrow\mathbb{\overline C}$ is continuous.
\QED

For pairs $(t,x)$ such that $p_t(x)=q_t(x)=0$ we shall see in the following lemma that 
for small $\epsilon>0$, $f_{t+\epsilon}$ is monotone on $(x,x+\epsilon)$ and the image 
of $(x,x+\epsilon)$ covers $\mathbb{\overline{R}}$ except for an interval of length 
$O(\epsilon)$. 

\begin{lem}\label{cont} 
Suppose that $y_1=0$ and that $\{(p_\epsilon,q_\epsilon):\epsilon\geq 0\}$ is a family of 
ordered pairs of linearly independent real polynomials of degree at most $d$ satisfying:
	
{\rm(i)} the coefficients of $p_\epsilon$ and $q_\epsilon$ are analytic in $\epsilon$,
	
{\rm(ii)} the roots of $w_\epsilon=W(p_\epsilon,q_\epsilon)$ are 
$y_1=0,y_2,\ldots ,y_{2d-3},\epsilon$,
	
{\rm(iii)} $p_0$ and $q_0$ have a common root at $x=0$, and
	
{\rm(iv)} $f_0(0)=0$.
	
Then for small $\epsilon >0$, $f_\epsilon(0)=O(\epsilon)$,
$f_\epsilon(\epsilon)=O(\epsilon)$ and $f_\epsilon$ has a pole in $(0,\epsilon)$. 
Moreover, $f_\epsilon:[0,\epsilon]\longrightarrow\mathbb{\overline{R}}$ is one-to-one and 
its image is $\mathbb{\overline{R}}\setminus (\alpha,\beta)$ where 
$(\alpha,\beta)$ is an open interval of size $O(\epsilon)$ with endpoints at 
$f_\epsilon(0)$ and $f_\epsilon(\epsilon)$.
\end{lem}

\textbf{Proof.} 
By our assumption, $p_0$ and $q_0$ have a common root at $x=0$, and $f_0(0)=0$, 
so there exist $k_1\neq 0$ and $k_2\neq 0$ such that $p_0=k_1x^2+O(x^3)$ and 
$q_0= k_2x+O(x^2)$. 
The requirement that $f_\epsilon$ has a fixed critical point at $x=0$ for all 
$\epsilon>0$ implies that 
$p_\epsilon(x)=k_1(x^2+a\epsilon x+ab\epsilon^2)+O(x^3,\epsilon x^2,\epsilon^2 x, \epsilon^3)$
and $q_\epsilon(x)=k_2(x+b\epsilon)+O(x^2,\epsilon x,\epsilon^2)$. 
Then $p'_\epsilon q_\epsilon-p_\epsilon q'_\epsilon=k_1k_2(x^2+2b\epsilon x)+
O(x^3, \epsilon x^2 ,\epsilon^2 x,\epsilon^3)$.

Since $f_\epsilon$ has a critical point at $x=\epsilon$, we should have $b=-\frac12$. 
Hence $q_\epsilon$ has a root at $x_\epsilon=\frac 1 2 \epsilon+O(\epsilon^2)$, 
between the critical points $0$ and $\epsilon$ of $f_\epsilon$, and 
$p_\epsilon(x_\epsilon)\neq 0$. 
Since $p_\epsilon(0)=O(\epsilon^2)$ and $p_\epsilon(\epsilon)=O(\epsilon^2)$, we have 
$f_\epsilon(0)=O(\epsilon)$ and $f_\epsilon(\epsilon)=O(\epsilon)$.

For small enough $\epsilon>0$, the interval $(0,\epsilon)$ does not contain any points of $D$.
If $f'_\epsilon<0$ on $(0,\epsilon)$ then when $x$ moves from $0$ to $\epsilon$, 
the value of $f_\epsilon(x)$ decreases from $O(\epsilon)$, tends to $-\infty$ as $x$ 
approaches $x_\epsilon$, and returns from $\infty$ to $O(\epsilon)$ at $x=\epsilon$. 
If $f'_\epsilon>0$ on $(0,\epsilon)$ then when $x$ moves from $0$ to $\epsilon$, the value of 
$f_\epsilon(x)$ increases from $O(\epsilon)$, tends to $\infty$ as $x$ approaches 
$x_\epsilon$, and returns from $-\infty$ to $O(\epsilon)$ at $x=\epsilon$.
\QED

Suppose $t\longrightarrow g_t=p_t^*/q_t^*$ is another path satisfying 
$L(p_t^*,q_t^*)=\chi(t)=L(p_t,q_t)$. 
There exists a real fractional linear transformation $h_t=(a_tx+b_t)/(c_tx+d_t)$ such that 
$g_t=h_t\circ f_t$. 
In particular, for each $t$, $f_t$ and $g_t$ have the same critical points. 
Observe that continuity of $p_t,q_t,p_t^*,q_t^*$ implies that the sign of $a_td_t-b_tc_t$ 
does not change as $t$ varies in $\mathcal{I}$. 
Post-composition with a real fractional linear transformation defines an equivalence 
relation on $\mathcal{{R}}_d$. This proves the following.

\begin{thm}\label{uniquec} With the above notation, let $f_{t_0}\in\mathcal{R}_d$ be a 
function whose Wronskian is equal to $\rho(t_0)$. 
Let $\mathcal{R}^o_{d}$ be the family of equivalence classes of functions in $\mathcal{R}_d$ 
of degree exactly $d$ and let $[f_{t_0}]\in\mathcal{R}_d^0$ be the equivalence class of 
$f_{t_0}$. 
There exists a unique map $\widehat{\eta}:\mathcal{I}\longrightarrow\mathcal{R}_d^0$ with 
$\widehat{\eta}(t_0)=[f_{t_0}]$ that projects to $\rho$.\end{thm}

The set $\mathcal{R}^0_{d}$ is in one-to-one correspondence with the set of nondegenerate 
real pencils. 
A pencil $L(p,q)\in G$ is nondegenerate if $\deg(p/q)=d$, in other words, 
if $\max(\deg p,\deg q)=d$ and $p$ and $q$ have no common zeros. 
Otherwise $L$ is said to be degenerate.

\begin{notn}\label{rmk} We shall identify $\mathbb{\overline{R}}$ with 
the unit circle $\mathbf{S}^1$ via a fractional linear transformation 
 $\vartheta:\mathbb{\overline{C}}\longrightarrow\mathbb{\overline{C}}$ 
mapping $\mathbb{\overline{R}}$ to $\mathbf{S}^1$, $0$ to $1$ and 
 preserving orientation. A real rational function $f$ will be identified 
 with the function $\widetilde{f}=\vartheta\circ f\circ\vartheta^{-1}$. 
 The critical points of $\widetilde{f}$ are the images of the critical 
 points of $f$ under $\vartheta$. In particular, all critical points of 
 $\widetilde{f}$ belong to the unit circle. Moreover, for any two points 
 $a$ and $b$ in $\mathbb{\overline{R}}$, $f(a)=f(b)$ if and only if 
 $\widetilde{f}(\vartheta(a))=\widetilde{f}(\vartheta(b))$. 

 Let $\mathcal{S}_d=\{\widetilde{f}=\vartheta\circ f\circ 
 \vartheta^{-1}:f\in\mathcal{R}_d\}$.
 We shall write $f$ instead of $\widetilde{f}$ and specify whether we are 
 considering $f$ as a function in $\mathcal{R}_d$ or $\mathcal{S}_d$. 
 With this identification, any construction or proof given for 
 $\mathcal{R}_d$ also applies to $\mathcal{S}_d$ and vice 		
versa. 
\end{notn}

 \begin{defn} Let $L=L(p,q)$ be a real pencil whose Wronskian has $2d-2$ 
 real roots counted with multiplicity and let $f=p/q$. The set 
$\Gamma=f^{-1}(\overline{\mathbb{R}})$ is
 called the net of $f$~\cite{EG1}.
Alternatively,  
let 
$\widetilde{f}=\vartheta\circ 
 f\circ\vartheta^{-1}$. Then $\Gamma=\widetilde{f}^{-1}(\mathbf{S}^1)$. 
\end{defn}

Observe that $\Gamma$ is independent of the choice of basis $(p,q)$, so 
we can also refer to $\Gamma$ as the net of $L$. 
Unless explicitly stated otherwise, we shall consider $\Gamma$ to be defined as 
$\widetilde{f}^{-1}(\mathbf{S}^1)$ rather than $f^{-1}(\overline{\mathbb{R}})$. 
Clearly, $\mathbf{S}^1\subset \Gamma$ and  $\Gamma$ is invariant under the 
reflection $s(z)=1/\overline z$ with respect to the unit circle.

The net $\Gamma$ defines a cell decomposition of $\mathbb{\overline{C}}$ 
as follows:

1) the vertices (0-cells) of the cell decomposition correspond to the zeros of 
$W(p,q)$,

2) the edges (1-cells) are the components of $\Gamma\setminus V$ where 
$V$ is the set of vertices,

3) the faces (2-cells) are the components of $\mathbb{C}\setminus\Gamma$.

This cell decomposition has the following properties. The closure of each 
cell is  homeomorphic to a closed ball of the corresponding dimension.
If $W(p,q)$ has distinct real roots, then each vertex of $\Gamma$ has degree 4.
If $W(p,q)$ has multiple roots we call $\Gamma$ a degenerate net. 

\begin{rem}\label{degendeg}  Let $L(p,q)$ be a real pencil whose net 
$\Gamma$ is degenerate and whose Wronskian has only simple or double 
roots, all real. A vertex of $\Gamma$ of degree 2 corresponds to a point 
$x$ such that $p(x)=q(x)=0$ while a vertex of degree 6 corresponds to a 
double root $x$ of $W(p,q)$ such that at least one of $p(x)$ and $q(x)$ 
is not 0. \end{rem}

\begin{defn} The edges of $\Gamma$ that lie inside the unit disc are called chords or interior edges of $\Gamma$. 
Since $\Gamma$ is symmetric with respect to $\mathbf{S}^1$, it suffices to consider its chords. Given a vertex $v$, the pair $\gamma=(\Gamma,v)$ is called the net of $f$ with respect to $v$, or simply the net of $f$ if $v$ is clear. The vertex $v$ is called the distinguished vertex of $\gamma$.
\end{defn}

The standard orientation of $\mathbf{S}^1$ induces a cyclic order $\prec$ on the vertices of $\Gamma$. We shall label these vertices $v_1,v_2,\ldots, v_{2d-2}$ so that $v_1\prec v_2\prec \ldots\prec v_{2d-2}\prec v_1$ and we shall set $v_{2d-1}=v_1$. For a net $\gamma=(\Gamma,v_1)$, we define a linear order on the vertices by $v_1<v_2<\ldots <v_{2d-2}$.

Two nets $(\Gamma,v_1)$ and $(\Gamma',v'_1)$ are said to be equivalent if there exists a homeomorphism $\mathbb{\overline{C}}\longrightarrow\mathbb{\overline{C}}$ mapping $\Gamma$ to $\Gamma'$, $v_1$ to $v'_1$, preserving orientation of both $\mathbb{\overline{C}}$ and $\mathbf{S}^1$ and commuting with the reflection $s$. By abuse of notation, we shall often refer to an equivalence class of nets by one of its representatives. 

\begin{defn} For a net $\gamma=(\Gamma,v)$, let $\hbox{\rm {Shift}}(\gamma)=(\Gamma,v')$ where $v'$ is the predecessor of $v$ under $\prec$.
\end{defn}

\begin{rem} Without loss of generality, we may assume that the vertices of $\Gamma$ 
are equally spaced on the circle. 
The net $\hbox{\rm Shift}(\gamma)$ is equivalent to the net $(\Gamma',v)$ 
where the ordered vertex set of $\Gamma'$ coincides with that of $\Gamma$ and the chords
of $\Gamma'$ are obtained by rotating the chords of $\Gamma$ by $\pi/(d-1)$ counterclockwise.
\end{rem}

\begin{thm}\label{unram} 
	The Wronski map $\mathcal{W}:G\longrightarrow\mathbb{P}^{2d-2}$ is unramified over the space $Q$ of all polynomials with real roots of multiplicity at most 2.
\end{thm}

\textbf{Proof.}  
Let $\pi_0\in Q$, $L(p_0,q_0)\in \mathcal{W}^{-1}(\pi_0)$ and $f_0=p_0/q_0$. 
The function $f_0$ can be chosen real~\cite{EG1}.
>From the definition of $Q$, the common roots $r_1,\dots,r_m$ of $p_0$ and $q_0$
(if any) must be simple. 
Suppose that $\mathcal{W}$ is ramified at $L(p_0,q_0)$, with ramification
index $\nu$. Let $\pi_j$ be a sequence of polynomials with distinct real
roots converging to $\pi_0$ as $j\to\infty$.
Since the polynomials with distinct real roots form an open
set in the space of all real polynomials, we can assume that 
$\mathcal{W}$ is not ramified over $\pi_j$ for all $j\ge 1$.
Then $\mathcal{W}^{-1}(\pi_j)$ contains $\nu$ distinct points $L(p_{j,k},q_{j,k})$
of the Grassmannian converging to $L(p_0,q_0)$ as $j\to\infty$.
Theorem \ref{eg} implies that all pencils $L(p_{j,k},q_{j,k})$ are real.
With the proper normalization, we can assume that $p_{j,k}$ and $q_{j,k}$ are
real polynomials, and that $(p_{j,k},q_{j,k})\to (p_0,q_0)$ as $j\to\infty$.
Then the rational functions $f_{j,k}=p_{j,k}/q_{j,k}$ are real, with
all real critical points, converge to $f_0$ as $j\to\infty$
uniformly on every compact set not containing the points $r_1,\dots,r_m$,
and their nets converge to the (degenerate) net of $f_0$.
Since the roots of $\pi_0$ are at most double, the net of $f_0$
may have vertices of degree either 6 or 2.
In both cases, the net of $f_{j,k}$ is uniquely determined by the net
of $f_0$ for large enough $j$, independent of $k$.
It follows from \cite{EG1} that, for large enough $j$,
pencils $L(p_{j,k},q_{j,k})$ do not depend on $k$.
Hence $\nu=1$ and $\mathcal{W}$ is not ramified over $\pi_0$. \QED

\section {  Properties of the nets of functions $f_t$.}
In this section, we study dependence on the parameter $t$
of the net $\gamma_t$ of the rational function $f_t$ defined in Section 2.
In particular, we show that the net $\gamma_t$ is shifted counterclockwise
when $\theta(t)$ passes its first vertex $y_1$.
We describe also degenerate nets corresponding to polynomials
$w_\theta$ with double roots.

\begin{defn} A $2\times(d-1)$ Young tableau is the distribution of the 
integers $1,2,\ldots,2d-2$ on a $2\times(d-1)$ rectangular array such 
that every integer is greater than the one above it and the one to its 
left, if any.
\end{defn}

Equivalence classes of nets with $2d-2$ vertices can be identified with 
$2\times(d-1)$ Young tableaux: 
to a net $\gamma=(\Gamma,v_1)$ with the ordered vertex set $\{v_1,v_2,\dots , v_{2d-2}\}$, 
we associate the tableau with an integer $i$ 
belonging to the first row if and only if $v_i$ is connected by a chord to a vertex 
$v_j$ with $j>i$. 
In particular, the number of equivalence classes of nets having $2d-2$ vertices is 
the $d$-th Catalan number $u_d=\frac{1}{d}{2d-2\choose d-1}$~(see \cite{St}). 

\begin{figure}[h]
	\includegraphics{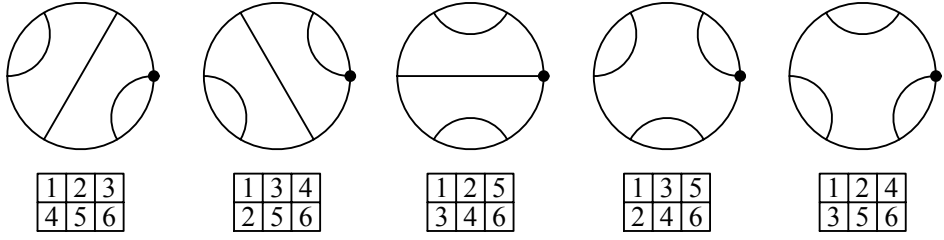}
	\caption{Nets and their Young tableaux for $d=4$.}
	\label{fig:d=4}
\end{figure}

\begin{ex} For $d=4$, there are $u_4=5$ equivalence classes of nets with $6$ vertices. 
These equivalence classes {\rm (}with the distinguished vertex $v_1=1${\rm )} and their respective 
Young tableaux are shown in {\rm Fig.~\ref{fig:d=4}.}\end{ex}

Let $D=\{y_1,y_2,\ldots,y_{2d-3}\}$ be an ordered set of distinct points in 
$\mathbf{S}^1$ all different from $1$. 
Let $\theta:[0,2d-2]\longrightarrow\mathbf{S}^1$ be the path given by $\theta(t)=e^{2\pi it}$.
By Lemma \ref{lifts} and Notation \ref{rmk}, there exists a family 
$\{f_t\in\mathcal{S}_d:t\in[0,2d-2]\}$ with the following properties:

i) for each $t$, the critical points of $f_t$ counted with multiplicity are 
$y_1,y_2,\ldots,y_{2d-3}$ and $\theta(t)$;

ii) given $(t_0,z_0)\in [0,2d-2]\times\overline{\mathbb{C}}$, $f_t(z_0)$ as a function of $t$ 
is continuous at $t_0$ whenever $\deg f_{t_0}=d$. 
In particular, for any $z_0\notin D$, $f_t(z_0)$ is a continuous function of $t$.

\bigskip

For $t\in[0,2d-2]$, let $\gamma_t=(\Gamma_t,y_1)$ be the net of $f_t$ with respect to $y_1$. 
Let $\theta^{-1}(D)=\{t_1, t_2,\ldots, t_{\kappa-1}\}$ with 
$0<t_1<t_2<\ldots<t_{\kappa-1}<2d-2$ and let $t_0=0$ and $t_\kappa=2d-2$. 
Here $\kappa-1=(2d-2)(2d-3)$ since $\theta$ makes $2d-2$ turns around $\mathbf{S}^1$ as $t$ 
goes from $0$ to $2d-2$. 
The sequence $\{s_1,s_2,\ldots,s_{2d-2}\}=\theta^{-1}(y_1)$ is a subsequence of $\{t_i\}$.
Let $s_0=0$ and $s_{2d-1}=2d-2$. 

\begin{lem}\label{samenet} Let $l\in\{0,1,\dots,2d-2\}$. 
For $t\in(s_l,s_{l+1})\setminus \theta^{-1}(D)$, the nets $\gamma_t$ are equivalent.
\end{lem}
\textbf{Proof.} Let $a$ and $b$ be two points in $(s_l,s_{l+1})\setminus \theta^{-1}(D)$ 
with $a<b$. 
If $\theta([a,b])\cap D=\emptyset$, then for each $z\in\mathbb{\overline{C}}$, $f_t(z)$ 
is a continuous function of $t$ on $[a,b]$ and, for all $t\in[a,b]$ the nets $\gamma_t$ 
are nondegenerate. 
Therefore $v_iv_j$ is a chord in $\gamma_{a}$ if and only if it is also a chord in 
$\gamma_{b}$, and hence the two nets are equivalent.

If $\theta([a,b])\cap D=\{y_i\}$, $i\neq 1$, then $\theta(a)=v_i$ and $y_i=v_{i+1}$ in 
the ordered vertex set of $\gamma_{a}$ while $y_i=v_{i}$ and $\theta(a)=v_{i+1}$ in that 
of $\gamma_{b}$. 
For each $z\in\mathbb{\overline{C}}\setminus\{y_i\}$, $f_t(z)$ is a continuous function 
of $t$ on $[a,b]$ and, for all $t\in[a,b]\setminus\theta^{-1}(y_i)$ the nets $\gamma_t$ 
are nondegenerate. 
Thus, for all $m,n\notin\{i,i+1\}$, $v_mv_n$ is a chord in $\gamma_{a}$ if and only if 
it is also a chord in $\gamma_{b}$. 
We have to consider the following two cases.

\begin{figure}[p]
	\includegraphics{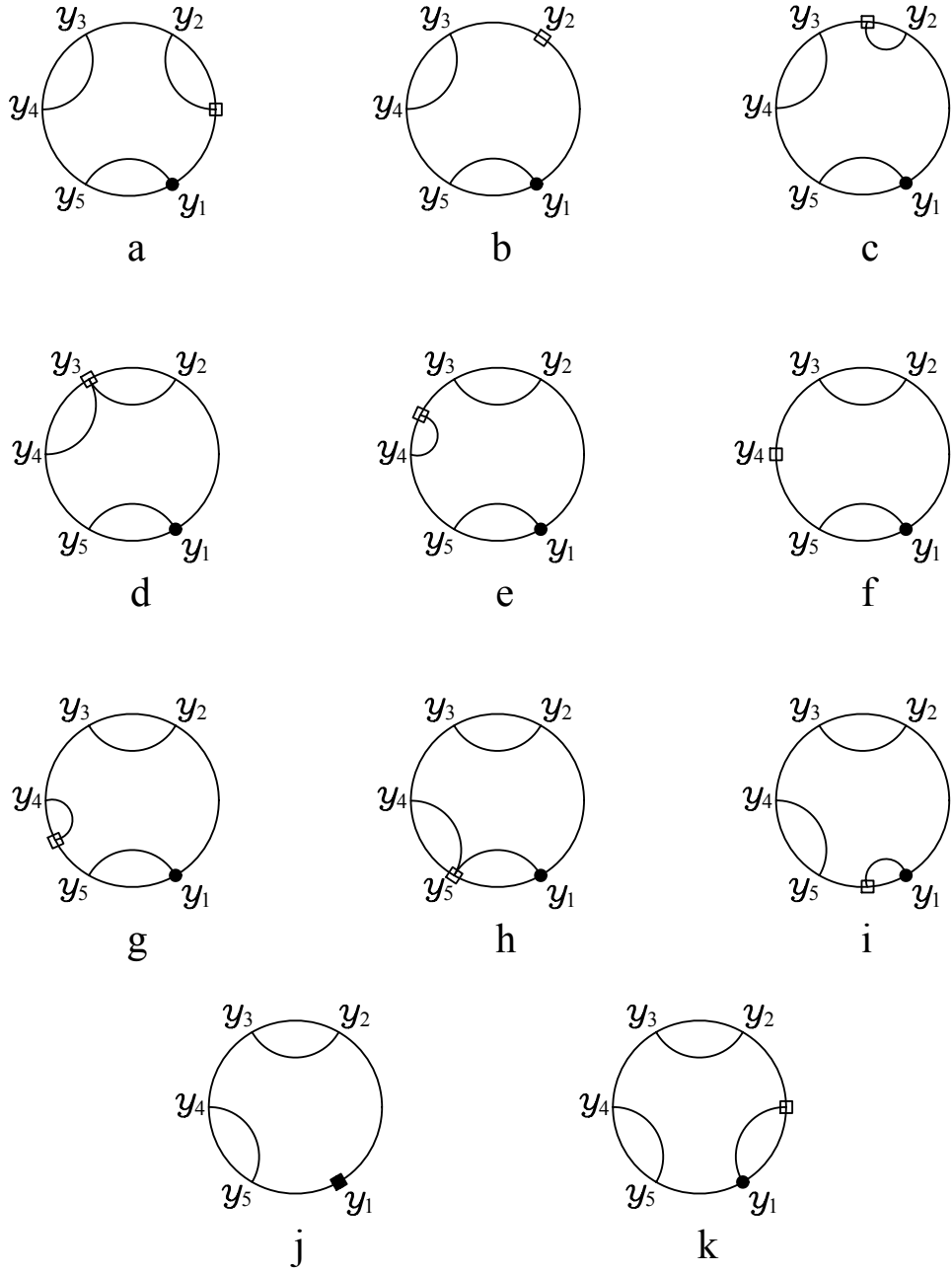}
	\caption{Dependence of the net $\gamma_t$ on $\theta=e^{2\pi it}$.}
	\label{fig:rotate}
\end{figure}

\textbf{Case 1:} $v_iv_{i+1}$ is a chord in $\gamma_{a}$. 
In this case, every other chord $v_mv_n$ in $\gamma_{a}$ satisfies $m,n\notin\{i,i+1\}$ 
and hence must also be a chord in $\gamma_{b}$. 
Therefore, $v_iv_{i+1}$ must be a chord in $\gamma_{b}$. 
An example of this case is given in Fig.~\ref{fig:rotate}a-c.

\textbf{Case 2:} $v_iv_{i+1}$ is not a chord in $\gamma_{a}$. 
Let $y_j$ be the endpoint of the chord from $v_i$ and $y_k$ the endpoint of that from 
$v_{i+1}$. 
Continuity of $f_t(z)$ on $[a,b]$ for any $z\neq y_i$ implies that $y_jy_k$ is not a chord 
in $\gamma_{b}$. 
Determining $\Gamma_{b}$ is thus reduced to finding a net with the vertex set 
$\{y_j\prec v_i\prec v_{i+1}\prec y_k\}$ such that $y_j$ and $y_k$ are not endpoints of 
the same chord. There are only two nets having four vertices, and exactly one of them 
does not have a chord joining $y_j$ to $y_k$. 
For an example, see Fig.~\ref{fig:rotate}c-e. \QED

\begin{rem}\label{degen} A similar argument can be used to determine the degenerate net 
$\Gamma_{t_i}$.
If, for $t$ close to $t_i$, there is a chord of $\Gamma_t$ connecting $\theta(t)$ and 
$\theta(t_i)$ {\rm (see, for example Fig.~\ref{fig:rotate}a,c),} then the vertex 
$\theta(t_i)$ is of degree 2 in $\Gamma_{t_i}$ {\rm(Fig.~\ref{fig:rotate}b).} 
Otherwise, $\theta(t_i)$ is of degree $6$ in $\Gamma_{t_i}$ {\rm(Fig. \ref{fig:rotate}d).}
\end{rem}

\begin{lem}\label{shift}  Let $[\gamma]_l=([\Gamma]_{l},y_1)$ denote the class 
$\{\gamma_t: t \in (s_{l},s_{l+1})\setminus\theta^{-1}(D)\}$. 
Then $[\gamma]_{l+1}=\hbox{\rm Shift}([\gamma]_l)$.
\end{lem}

\textbf{Proof.} Choose representatives $\gamma_{t'}$ and $\gamma_{t''}$ of $[\gamma]_l$ and 
$[\gamma]_{l+1}$ repectively with $y_{2d-2}\prec\theta(t')\prec y_1\prec\theta(t'')\prec y_2$.
As $t$ goes from $t'$ to $t''$, $\theta(t)$ crosses the distinguished vertex $y_1$, so 
Lemma \ref{samenet} cannot be applied to $\gamma_{t'}$ and $\gamma_{t''}$. 
Instead we shall consider $\gamma'=(\Gamma_{t'},y_2)$ and $\gamma''=(\Gamma_{t''},y_2)$ 
which are related to $\gamma_{t'}$ and $\gamma_{t''}$ by $\gamma_{t'}=\hbox{Shift}(\gamma')$ 
and $\gamma_{t''}=\hbox{Shift}^2(\gamma'')$. 
Since $\theta(t)$ does not cross $y_2$ as $t$ goes from $t'$ to $t''$, we can apply Lemma 
\ref{samenet} to $\gamma'$ and $\gamma''$ to get that $\gamma'=\gamma''$. 
Thus $\gamma_{t''}=\hbox{Shift}^2(\gamma'')=\hbox{Shift}^2(\gamma')=\hbox{Shift}(\gamma_{t'})$.
\QED

\begin{ex} 
The nets a, c, e, g in {\rm Fig.~\ref{fig:rotate}} all belong to $[\gamma]_0$ while 
the net k belongs to $[\gamma]_{1}$ and is the shift of $[\gamma]_0$.
\end{ex}

\begin{cor}\label{equalnetst} If $t'\in(t_0,t_1)$ and $t''\in(t_{\kappa-1},t_\kappa)$ 
then $\gamma_{t'}=\gamma_{t''}$. 
Equivalently, $[\gamma]_0=[\gamma]_{2d-2}$.\end{cor}

\textbf{Proof.} As $t$ goes from $t'$ to $t''$, $\theta(t)$ makes $2d-2$ turns around 
the circle. 
By Lemma \ref{shift}, each complete turn of $\theta$ corresponds to the shift operator 
(i.e., rotation by $\pi/(d-1)$) applied to the net. 
Applying the shift operator $2d-2$ times results in a complete rotation of the net, 
resulting in the original net.\QED

\begin{rem} It could happen that $[\gamma]_0=[\gamma]_l$ for some positive integer $l<2d-2$. 
This occurs when $[\gamma]_0$ has rotational symmetry of order $n$ for some integer $n>1$, 
in which case $l=(2d-2)/n$, a factor of $2d-2$. 
For example, the net $[\gamma]_0$ in {\rm Fig.~\ref{fig:rotate}a} has rotational symmetry of 
order $3$, so $[\gamma]_0=[\gamma]_{2}$.
\end{rem}

\begin{ex}
{\rm Fig.~\ref{fig:rotate}} describes how a given net is modified as $\theta(t)$ makes a 
full turn around the circle.  
As $\theta(t)$ moves in $(y_1,y_2)$ {\rm(Fig.~\ref{fig:rotate}a),} the net remains unchanged. 
When $\theta(t)$ reaches $y_2$ {\rm(Fig.~\ref{fig:rotate}b),} the net becomes degenerate 
with $\deg(y_2)=2$. 
As $\theta(t)$ moves away from $y_2$ {\rm(Fig.~\ref{fig:rotate}c),} we recover the original 
net until $\theta(t)$ reaches $y_3$ {\rm(Fig.~\ref{fig:rotate}d).}
At this point the net becomes degenerate again with $\deg(y_3)=6$. 
As $\theta(t)$ moves away from $y_3$ {\rm(Fig.~\ref{fig:rotate}e),} we recover the original 
net and the process continues in a similar manner until $\theta(t)$ crosses $y_1$. 
The net obtained at this step {\rm(Fig.~\ref{fig:rotate}k)} is the shift of the net in 
{\rm Fig.~\ref{fig:rotate}a.} 
\end{ex}

\section{Lower bounds for $k=1$.}

In this section, we derive lower bounds on the
number of real solutions to Problem \ref{rrf} in the special case
when the interval $(r,s)$, or its complement, contains only one fixed 
critical point of $f$.
 To simplify notation, we shall assume that $y_1=\infty$ is the only 
fixed vertex not contained in the interval $(r,s)$. Throughout this section, $f$ is a real rational function of degree $d$ and its net is given by $\Gamma=f^{-1}({\mathbb{\overline R}})$.
 
\begin{thm}\label{1pt} 
Let $r<y_2<y_3<\dots < y_{2d-3}<s$ be points in $\mathbb{R}$ and let $y_1=\infty$. Then there are at least $u_d-2u_{d-1}$ classes of real rational functions $f$ of degree $d$ having critical points at $y_j$, $1\leq j\leq 2d-3$ and satisfying $f(r)=f(s)$.
\end{thm}

\noindent\textbf{Proof} 
Let $\mathcal{C}$ be an equivalence class of the nets with $2d-2$ vertices and no interior 
edges connecting the distinguished vertex $y_1$ to any of its two neighboring vertices. 
The number of such classes is $u_d-2u_{d-1}$, since there are $u_{d-1}$ classes of nets 
having two given neighboring vertices connected by an edge. 
It remains to show that for each such class $\mathcal{C}$ there exists $y\in\mathbb{R}$ and 
a real rational function $f_y$ having the net $\gamma_y\in \mathcal{C}$ with the vertex set 
$\{y_1,y_2,\ldots,y_{2d-3},y\}$, satisfying $f_y(r)=f_y(s)$.

Let $y\in (s,\infty)$. There exists a unique class $\mathcal{F}_{y}$ of 
real rational functions with the net belonging to $\mathcal{C}$ and 
 critical points $y_1,y_2,\ldots, y_{2d-3},y$~\cite{EG1}. Choose $f_{y}\in 
 \mathcal{F}_{y}$ so that $f_{y}$ has a double pole at $\infty$ and 
 $f_{y}(x)<0$ for large $|x|$. 
Let $\rho:\mathbb{R}\longrightarrow Q$ be given by 
 $\rho(y)=(x-y_2)(x-y_3)\dots (x-y_{2d-3})(x-y)$. 
Note that $\rho(y)$ is the Wronskian of $f_{y}$. 
By Lemma \ref{lifts}, $\rho$ can be lifted 
to a path $\eta:\mathbb{R}\longrightarrow \mathcal{R}_d$ with 
$\eta(y)=f_{y}$ satisfying properties (i)-(iii) of Lemma \ref{lifts}. 
By Lemma \ref{samenet}, $\gamma_y\in\mathcal{C}$ for all 
 $y\in\mathbb{R}-\{y_j\}_{j=2}^{2d-3}$.

 The map $\phi_x:\mathbb{R}\longrightarrow \mathbb{\overline{R}}$ given by 
 $\phi_x(y)=f_y(x)$ is continuous for all $x$ in 
 $\mathbb{R}-\{y_j\}_{j=2}^{2d-3}$. Continuity of $\phi_x(y)$ on 
 $(s,\infty)$ implies that for all $y\in (s,\infty)$, $f_y(x)<0$ for 
 large $|x|$. The following lemma completes the proof.

\begin{lem} 
	There exists $y\in\mathbb{R}$ such that $f_y(r)=f_y(s)$.
\end{lem}

 \textbf{Proof.} Assume that $f_y(r)\neq f_y(s)$ for all $y\in [r,s]$. If 
 $y\in [r,s]$ then  $f_y$ cannot have a pole at $r$ nor at $s$  since each 
 of $[s,\infty)$ and $(-\infty,r]$ belongs to the boundary of a 
 face of $\Gamma_y$ and $f_y$ has a pole at $\infty$. Since both $\phi_r$ 
 and $\phi_s$ are continuous, this implies that either 
 $-\infty<f_y(r)<f_y(s)<\infty$ for all $y\in [r,s]$ or 
 $-\infty<f_y(s)<f_y(r)<\infty$ for all $y\in [r,s]$. Without loss of 
 generality, we may assume that $-\infty<f_y(r)<f_y(s)<\infty$ for all 
 $y\in [r,s]$. In particular, $-\infty<f_s(r)<f_s(s)<\infty$.

  For $y\in (s,\infty)$, $f_y(y)$ is finite since $f_y$ has a pole at 
 $\infty$ and $(y,\infty)$ belongs to the boundary of a face of 
 $\Gamma_y$. Moreover, $f_y(s)$ cannot exceed $f_y(y)$ without first 
 decreasing to $-\infty$ since $\phi_s$ is continuous and $f_y$ has a 
 local maximum at $y$ and no critical points between $s$ and $y$. On the 
 other hand, since $\Gamma_y$ has no interior edges connecting $\infty$ to 
 $y$ for all $y\in [s,\infty)$, $\lim_{y\rightarrow\infty} f_y(y)=-\infty$ 
 and it follows that $f_y(s)$ must also decrease to $-\infty$ on a 
 subinterval of $(s,\infty)$. In particular, there exists $y\in 
 (s,\infty)$ such that $f_y(r)=f_y(s)$.\QED

\section{Upper and lower bounds for the arc length of $f([r,s])$}
In this section, we assume that $f$ maps $\mathbf{S}^1$ to $\mathbf{S}^1$. 
We derive lower and upper bounds
for the arc length of $f([r,s])$, in terms of the net of $f$ and the 
position of the points $r$ and $s$ relative to the vertices of the net.

 Consider $\mathbf {S}^1$ with the standard orientation. Given any two 
 points $r$ and $s$ in $\mathbf{S}^1$, we define $[r,s]$ and $(r,s)$ 
 respectively to be the closed and open positively oriented arcs in 
 $\mathbf{S}^1$ starting at $r$ and ending at $s$. The standard 
 orientation on $\mathbf{S}^1$ induces an order on $[r,s]$ given by 
 $a\leq b$ if and only if $a\in [r,b]$. 

 \begin{defn}\label{YST} Let $\Gamma$ be a net such that
all its vertices $y_1<y_2<\ldots<y_n$ inside $(r,s)$ are simple. 
The Young tableau of $\Gamma$ corresponding to $(r,s)$, denoted by $Y_{\Gamma}(r,s)$, 
has $2$ rows and $n$ entries defined as follows. 
The integer $i$ is placed in the first row of $Y_{\Gamma}(r,s)$
if and only if $y_j$ is not connected by a chord of $\Gamma$ to $y_i$, $i>j$. 
\end{defn}

\begin{figure}[h]\begin{center}
	\includegraphics{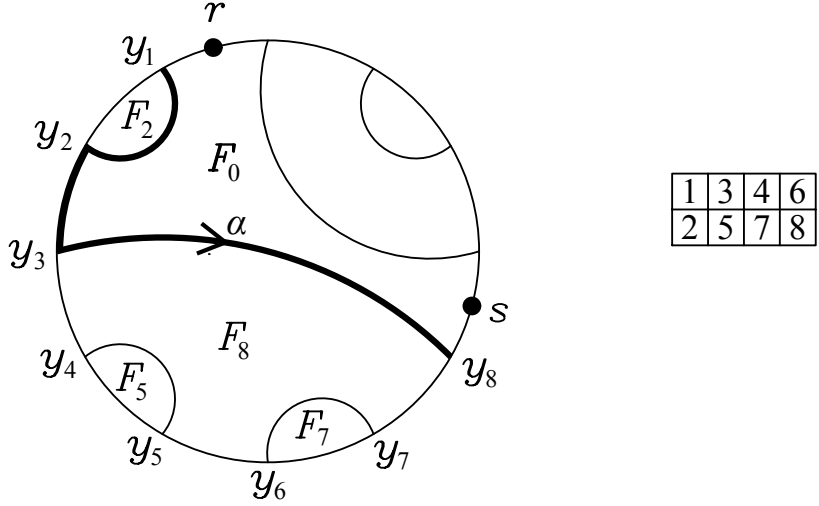}
 	\caption{A net $\Gamma$ with $m=0$ and the corresponding Young 
tableau $Y_\Gamma(r,s)$.}
	\label{speccase}\end{center}
\end{figure}

\begin{figure}[h]\begin{center}
	\includegraphics{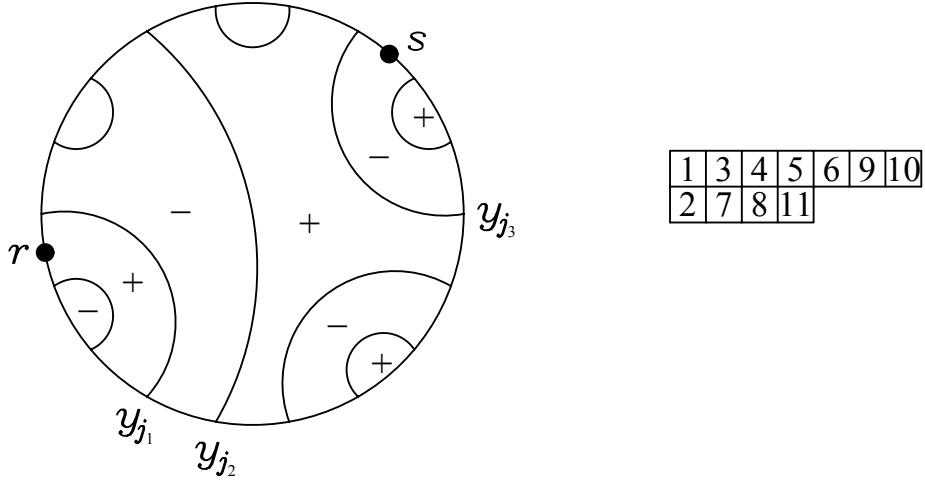}
 	\caption{A net $\Gamma$ with $m=3$ and the corresponding 
Young tableau $Y_{\Gamma}(r,s)$.}
	\label{gencase}\end{center}
\end{figure}

Note that $Y_{\Gamma}(r,s)$ is part of the Young tableau of 
$(\Gamma,y_1)$. See examples in Fig.~\ref{speccase} and Fig.~\ref{gencase}.

\begin{defn}Given a rational function $f$ and an oriented segment or 
a closed loop $c$ in $\mathbb {\overline{C}}$ with  $f(c)\subset 
\mathbf{S}^1$, let $L(c)=L_{f}(c)$ be the argument  increment of $f$ on 
$c$.
\end{defn}

Let $\Gamma$ be a net with simple vertices $y_1<\dots<y_n$ inside $(r,s)$. 
Let $f$ be a rational function with the net $\Gamma$. 
Assume that $f$ is orientation preserving on $(r,y_1)$.
Let $Y=Y_\Gamma(r,s)$.

\begin{defn}\label{EO}Let $E$ and $O$ be the 
numbers of even and odd entries in the second row of  $Y$. 
Let $m$ be the number of vertices of $\Gamma$ inside $(r,s)$ connected by a chord 
to vertices outside $(r,s)$.\end{defn}

\begin{lem} The number $m=n-2(O+E)$ is the difference between the length
of the first row of $Y$ and the length of its second row.\end{lem}

Consider first the case $m=0$.

\begin{lem}\label{(*)}The following are equivalent:
\begin{enumerate}[$(i)$]
\item $m=0$
\item $Y$ is rectangular
\item no vertex in $(r,s)$ is connected by a chord of $\Gamma$ to a vertex outside $(r,s)$.
\end{enumerate}
\end{lem}
 \begin{defn} Assuming $m=0$, let $V=\{y_1,\ldots,y_n\}$ be the set of vertices of $\Gamma$
inside $(r,s)$. 
Let $C$ be the set of chords of $\Gamma$ having both endpoints in $V$ and let 
$\mathcal{F}=\mathcal{F}(r,s)$ be the set of faces $F$ of $\Gamma$ 
satisfying $\partial F\cap \mathbf{S}^1\subset [r,s]$. A chord in $C$ 
 will be denoted by $y_iy_j$ where $y_i$ and $y_j$ are its endpoints, with 
 $i<j$. Given a chord $y_iy_j\in C$, let $F_j$ be the face in 
 $\mathcal{F}$ having $y_j$ as its largest vertex.  Let $F_0$ be the face 
 of $\Gamma$ whose boundary contains $[r,y_1]$.\end{defn}
 There is a one-to-one correspondence between the sets $C$ and 
 $\mathcal{F}$. Condition (iii) of Lemma~\ref{(*)} implies that $\partial 
 F_0$ must also contain $[y_n,s]$.
 \begin{ex}  For the net $\Gamma$ and arc $(r,s)$ in {\rm Fig.~\ref{speccase},}
 $n=8$, $E=2$, $O=2$, $m=0$,  $C=\{y_1y_2, y_3y_8, y_4y_5, 
 y_6y_7\}$ and $\mathcal{F}=\{F_2,F_5,F_7,F_8\}$.\end{ex}
   \begin{defn} Assuming $m=0$, a chord $y_iy_j$ in $C$ is even (resp., odd) if $j$ is 
even (resp., odd). \end{defn}
  
  By the construction of $Y_{\Gamma}(r,s)$ we have the following.
 
 \begin{lem} The numbers $E$ and $O$ in Definition \ref{EO} are, respectively, 
the numbers of even and odd chords in $C$.\end{lem}

\begin{defn}  A face $F$ of $\Gamma$ is positive if $f$ is orientation 
preserving on  $\partial F$, otherwise $F$ is negative. Let $P$ and $N$ be 
the numbers of  positive and negative faces in $\mathcal{F}$. The face 
$F_0$ is always positive since $f$ is orientation preserving on $(r,y_1)$.\end{defn}

Since $f$ maps the boundary of each face $F\in \mathcal{F}$ bijectively  
onto $\mathbf{S}^1$, it follows that $L(\partial F)=2\pi$ if  $F$ is 
positive and $L(\partial F)=-2\pi$ if $F$ is negative.
  
\begin{lem} Assuming $m=0$, a chord $y_iy_j$ in $C$ is odd if and only if $F_j$ 
is positive. In particular, $O=P$ and $E=N$.
\end{lem}
  
\textbf{Proof.}
For $y_iy_j\in C$, the arc $[y_{j-1},y_j]$ of $\mathbf{S}^1$ belongs to $\partial F_j$, 
and its orientation induced by $F_j$ agrees with the standard orientation on $\mathbf{S}^1$. 
Thus, since $f$ preserves orientation on $[r,y_1]$ and has only simple critical points in 
$(r,s)$, it preserves orientation on $[y_{j-1},y_j]$ if and only if $j$ is odd. \QED
  
\begin{lem}\label{posface} Let $m=0$. Then $L([r,s])\in(2\pi(O-E),2\pi(O-E+1))$.
\end{lem}

\textbf{Proof.} Let $\Omega$ be 
 the region bounded by the arc $[y_1,y_n]\subset\mathbf{S}^1$ and 
 $\alpha$, the positively oriented curve on $\partial F_0$ from $y_1$ to 
 $y_n$ (see Fig.~\ref{speccase}). Since  $\Omega\setminus\Gamma$ is the 
union of the faces in $\mathcal{F}$, $$L(\partial\Omega)=\sum_{F\in 
\mathcal{F}}L(\partial F)=2\pi(P-N)=2\pi(O-E).$$
 On the other hand, $L([r,y_1])+L(\alpha)+L([y_n,s])\in(0,2\pi)$ since 
 $[r,y_1]\cup\alpha\cup[y_n,s]$ is part of the boundary of the positive face $F_0$ 
 having positive orientation. Therefore,
 $$\begin{array}{lcr}
 L([r,s])&=& L([r,y_1])+L([y_1,y_n])+L([y_n,s])\\
 &=& L([r,y_1])+L(\alpha)+ L(\partial\Omega)+L([y_n,s])\\
&\in&(2\pi(O-E),2\pi(O-E+1)).
\end{array}$$
\QED

\begin{rem}\label{negface} 
Assume $m=0$. If $f$ is orientation reversing on $(r,y_1)$ then  $g=1/f$ is 
orientation preserving 
on $(r,y_1)$ and $L_g([r,s])=-L_f([r,s])$. Replacing $f$ by $g$ we obtain 
from Lemma 
\ref{posface}  
$L_f([r,s])\in(2\pi(E-O-1),2\pi(E-O))$.\end{rem}

Consider now the case $m>0$.   

\medskip
 Let $y_{j_1},y_{j_2},\ldots,y_{j_m}$ be the subsequence of 
 $y_1,y_2,\ldots, y_n$ consisting of the vertices joined by a chord 
of $\Gamma$ to vertices outside $(r,s)$.
  Letting $j_0=0$, $y_{j_0}=r$, $j_{m+1}=n+1$ and $y_{j_{m+1}}=s$, we can 
represent 
 $[r,s]$ as the union $\bigcup_{i=0}^m[y_{j_i},y_{j_{i+1}}]$. 
Let $i\in\{0,1,\ldots,m\}$ and let $E_i$ and $O_i$ be the 
 numbers of even and odd entries in the second row of $Y$ that are between 
 $j_i$ and $j_{i+1}$.
 \begin{ex} For the net $\Gamma$ and arc $(r,s)$ in {\rm Fig.~\ref{gencase},}
 $E=O=2$, $m=3$, $E_0=1$, $O_0=0$, $E_1=O_1=0$, $E_2=O_2=1$, $E_3=0$, 
 $O_3=1$.\end{ex}

For the arc $(y_0,y_{j_1})$, $f$ is orientation preserving on $(y_0,y_1)$, so by 
Lemma \ref{posface},   
$L([y_0,y_{j_1}])\in(2\pi(O_0-E_0),2\pi(O_0-E_0+1))$. For 
$(y_{j_1}, y_{j_2})$,  $f$ is orientation reversing on $(y_{j_1}, y_{j_1+1})$. Since 
$j_1$ is odd, the vertices that are even with respect to $(y_{j_1}, 
y_{j_2})$ belong to $O_1$ while the vertices that are odd with respect to 
$(y_{j_1}, y_{j_2})$ belong to $E_1$. Accordingly, by Remark 
\ref{negface},  
$L_f([y_{j_1},y_{j_2}])\in(2\pi(O_1-E_1-1),2\pi(O_1-E_1))$. 
Applying the same arguments to $(y_{j_i},y_{j_{i+1}})$ we get the 
following.

  \begin{lem}  Let $i\in\{0,1,\ldots,m\}$. If $i$ is even then  $L([j_i,j_{i+1}])\in (2\pi(O_i-E_i),2\pi(O_i-E_i+1))$. If $i$ is odd then  $L([j_i,j_{i+1}])\in (2\pi(O_i-E_i-1),2\pi(O_i-E_i))$.\end{lem}
  
Summing up over $i=0, \ldots ,m $ we obtain the following.
  
\begin{thm}\label{igammat} Let $I=I_{\Gamma}(r,s)=$
$$\big([(n+1)\hskip-0.1in\mod 2 -n+4O-1]\pi,\,[(n+1)\hskip-0.1in\mod 2 +n-4E+1]\pi\big).$$
Then $L([r,s])\in I$.\end{thm}
  
Now we consider the case when $W(p,q)$ has a double root $y$ and the net 
$\Gamma$ is degenerate. If $y\notin (r,s)$ this does not affect 
$L([r,s])$. If $y\in (r,s)$ then in labelling the vertices $y_1,\ldots, 
y_n $ in $(r,s)$, we 
assign to $y$ two consecutive indices and define $Y_\Gamma(r,s)$ to be $Y_{\Gamma'}(r,s)$ 
where $\Gamma'$ is a nondegenerate net with the vertices and chords close
to those of $\Gamma$. With this agreement all the 
arguments above can be applied to $\Gamma'$ and the statement of 
\ref{igammat} remains true except for the case when the degree of $y$ is 6 and 
both chords of $\Gamma$ with the ends at $y$ have other ends outside $(r,s)$. 
This means that $y$ represents a segment $(y_{j_i},y_{j_{i+1}})$ of length zero.

\begin{thm} \label{igammatd} Let $y$ be a double root of $W(p,q)$ inside $(r,s)$ which is a 
vertex of $\Gamma$ of degree 6 such that both chords of $\Gamma$ with the ends at $y$ 
have other ends outside $(r,s)$. 
Let $I=I_{\Gamma}(r,s)$ be
$$\big([(n+1)\hskip-0.1in\mod 2-n+4O-1]\pi,\,[(n+1)\hskip-0.1in\mod 2+n-4E-1]\pi\big)$$ 
if the number of vertices between $r$ and $y$ is odd, and 
$$\big([(n+1)\hskip-0.1in\mod 2-n+4O+1]\pi,\,[(n+1)\hskip-0.1in\mod 2+n-4E+1]\pi\big)$$ 
otherwise. Then $L([r,s])\in I$. 
\end{thm}

We can similarly compute $I_{\Gamma}(s,r)$ and use it to improve $I_{\Gamma}(r,s)$ in some cases.

Let $f$ be a function in $S_d$ all of whose critical points are simple except for possibly one double critical point. Let the critical points of $f$ be $y_1, \dots, y_{2d-2}$ (counted with multiplicity) with $y_1\prec r\prec y_2\preceq y_3\preceq \dots \preceq y_{k+1}\prec s\prec y_{k+2}\preceq\dots \preceq y_{2d-2}\preceq y_1$. Assume $f$ is orientation preserving on $(y_1,y_2)$. Let $(\Gamma, y_1)$ be the net of $f$ with respect to $y_1$ and let $Y$ be the corresponding Young tableau. Let $\mathcal{F}$ be the set of all faces of $\Gamma$. 
Let $E$ and $O$ be the numbers of even and odd entries in the second row of $Y$ respectively. Then $$L([r,s])+L([s,r])=\sum_{F\in\mathcal{F}}\partial F=2\pi(E-O-1).$$  
So if $I_{\Gamma}(r,s)=(a, b)$,  $I_{\Gamma}(s,r)=(a', b')$ and $L([r,s])+L([s,r])=c$, then $I_{\Gamma}(r,s)$ can be improved to $(a,b)\cap (c-b', c-a')$.
 
\section{Lower bounds in the general case.}

In this section, we apply the results of Section 5
to derive lower bounds on the number of real solutions to Problem 1.5.
Without loss of generality, we may set $r=-1$ and $s=1$. 
Let $D=\{y_1,y_2,\ldots, y_{2d-3}\}$ be a set of distinct points in $\mathbf{S}^1$ all 
different from $\pm 1$.  
Assume that $y_1$ does not belong to the arc $[r,s]$. 
Let $f_0$ be a function in $\mathcal{S}_d$ having critical set $D\cup\{1\}$. 
Let $\theta:\mathbb{R}\longrightarrow\mathbf{S}^1$ be the path given by 
$\theta(t)=e^{2\pi it}$. 
By Lemma \ref{lifts} and Notation \ref{rmk}, there exists a path 
$\eta:\mathbb{R}\longrightarrow \mathcal{S}_d$ given by $\eta(t)=f_t$, 
where $f_t=p_t/q_t$ is a rational function having critical points 
$y_1,y_2,\ldots, y_{2d-3},\theta(t)$ counted with multiplicity, 
satisfying properties (i)-(iii) of Lemma \ref{lifts}. 

\begin{lem} The map $\mathbb{R}\longrightarrow\mathbb{R}$ given by $t\mapsto L_{f_t}([r,s])$ 
is continuous everywhere except on a subset $\mathcal{D}$ of $\theta^{-1}(D)$ given by 
$\mathcal{D}=\{t\in\theta^{-1}(D):p_t(x)=q_t(x)=0 \hbox{ for some } x\in D\}$. 
At any point of discontinuity $t^*$, $\lim\limits_{t\to t^*}L_{f_t}([r,s])$ 
exists and $|L_{f_{t^*}}([r,s])- \lim\limits_{t\to t^*}L_{f_t}([r,s])|=2\pi$.\end{lem}

\textbf{Proof.} The first statement follows from Lemma \ref{lifts}. 
The second statement follows from Lemmas \ref{lifts} and \ref{cont} and the fact that 
$f_{t^*}(x)$ is continuous at $\theta(t^*)$ as a function of $x$.\QED

Let $L:\mathbb{R}\longrightarrow\mathbb{R}$ be given by 
$$
L(t) = \left\{ \begin{array}{ll}             
	L_{f_t}([r,s]) & \quad\hbox{if } t\notin\mathcal{D}\\
	\lim\limits_{x\rightarrow t}L_{f_{x}}([r,s]) & \quad\hbox{if }t\in\mathcal{D}\\
		 \end{array}       
					\right. $$

\begin{cor}\label{eqrs} The map $L$ is continuous on $\mathbb{R}$. 
For any $t\in\mathbb{R}$, $L(t)\in 2\pi\mathbb{Z}$ if and only if $f_t(r)=f_t(s)$.\end{cor}

For $t\in\mathbb{R}$, let $\gamma_t=(\Gamma_t,y_1)$ be the net corresponding to $f_t$. 
The net $\Gamma_t$ is nondegenerate if and only if $\theta(t)\notin D$. 
If $t\in\mathcal{D}$ then $\theta(t)$ is a vertex of degree 2 in $\Gamma_t$, and if 
$t\in\theta^{-1}(D)\setminus\mathcal{D}$ then $\theta(t)$ is of degree 6 in $\Gamma_t$. 
This follows from Lemma \ref{lifts} and Remark \ref{degendeg}.

\begin{lem}\label{periodic} The map $t\mapsto\gamma_t$ is periodic. 
The period $T$ is a factor of $2d-2$.
It is equal to the number of distinct nondegenerate nets in the image.
\end{lem}

\textbf{Proof.} This follows from Corollaries \ref{samenet} and \ref{equalnetst} and 
Remark \ref{degen}. \QED

For $n\in\mathbb{Z}$, let $V_n=(n,n+\frac12)$ and $W_n=(n+\frac12,n+1)$ so that 
$\theta(V_n)=(s,r)$ and $\theta(W_n)=(r,s)$. 
For $t\in\mathbb{R}$, let $2\pi l_1(t)$ and $2\pi u_1(t)$ be the lower and upper endpoints of the 
interval $I_{\Gamma_t}(r,s)$ defined in Theorems \ref{igammat} and \ref{igammatd}, and let $2\pi l_2(t)$ and $2\pi u_2(t)$ be the lower and upper endpoints of the interval $I_{\Gamma_t}(s,r)$ defined similarly. Let $2\pi c=L([r,s])+L([s,r])$. Then $L(t)/(2\pi)\in (l_1(t),u_1(t))\cap(c-u_2(t),c-l_2(t))$.

\begin{lem} The functions $l_1(t)$ and $u_1(t)$ are constant on each set $V_n$.\end{lem}

\textbf{Proof.} The interval $I_{\Gamma_t}(r,s)$ depends only on the Young tableau of $\Gamma_t$ 
corresponding to $(r,s)$ which is independent of any changes outside $(r,s)$.\QED

\begin{lem}\label{L(t)} The functions $l_1(t)$ and $u_1(t)$ are constant on each set 
$W_n\setminus(\theta^{-1}(D)\setminus\mathcal{D})$ and assume only finitely many values on each set $W_n$.\end{lem}

\textbf{Proof.} The arc $(r,s)$ does not contain the distinguished vertex $y_1$, so by 
Corollary \ref{samenet}, the nets $\gamma_t$, $t\in W_n\setminus\theta^{-1}(D)$ are equivalent.
Since $\theta(t)\in(r,s)\setminus D$ for all $t\in W_n\setminus\theta^{-1}(D)$, 
Definition \ref{YST} implies that
the Young tableaux $Y_{\Gamma_t}(r,s)$ are identical for all 
$t\in W_n\setminus\theta^{-1}(D)$. For $t\in W_n\cap\mathcal{D}$, the Young tableaux 
$Y_{\Gamma_t}(r,s)$ are identical to those for $t\in W_n\setminus\theta^{-1}(D)$,
as defined after Theorem \ref{igammat}. This along with the fact that $W_n\cap(\theta^{-1}(D)\setminus\mathcal{D})$ is finite proves the second part of the statement.\QED

If we let $\gamma_t=(\Gamma_t,y)$ for some fixed vertex $y\in (r,s)$ we get the following.

\begin{cor} The functions $l_2(t)$ and $u_2(t)$ are constant on each set $W_n$ and assume only finitely many values on each set $V_n$.\end{cor}

Let $E=E(f_0)$ be the set of equivalence classes of functions $f_t$, $t\in\mathbb{R}$. The map $t\longrightarrow [f_t]$ is periodic with period $T$. The value of $T$ depends on the choice of $f_0$, but it is always a factor of $2d-2$. Therefore, instead of calculating $T$ for each $f_0$, we shall consider the interval $[0,2d-2]$ and deal with the issue of having counted some equivalence classes more than once at the end of this section.

Let $\mathcal{U}$ and $\mathcal{L}$ be two integer-valued functions defined  by 
$$\begin{array}{rcl}
\mathcal{U}(n)&=&\inf\left\{u_1(t),c-l_2(t):t\in V_n\right\}\\
\mathcal{U}(n+\frac12)&=&\inf\left\{u_1(t),c-l_2(t):t\in W_n\right\}\\
\mathcal{L}(n)&=&\sup\left\{l_1(t),c-u_2(t):t\in V_n\right\}\\
\mathcal{L}(n+\frac12)&=&\sup\left\{l_1(t),c-u_2(t):t\in W_n\right\}\\
\end{array}$$

The functions $\mathcal{U}$ and $\mathcal{L}$ can be easily computed since $l_1, u_1, l_2, u_2$ assume only finitely many values on each $V_n$ and $W_n$. In addition, continuity of $L(t)$ implies  the existence of $t_n\in V_n$ and $t_{n+1/2}\in W_n$ such that $\frac{1}{2\pi}L(t_n)\in(\mathcal{L}(n),\mathcal{U}(n))$ and $\frac{1}{2\pi}L(t_{n+1/2})\in(\mathcal{L}(n+1/2),\mathcal{U}(n+1/2))$.
Let $S=\{k/2:0\leq k<4d-4\}$. Let $\prec$ be the cyclic order on $S$ given by $0\prec 1/2\prec 1\prec \dots\prec (4d-5)/2\prec 0$. For $i,j\in S$, let $(i,j)=\{k\in S|i\prec k\prec j\}$. 

\begin{defn} A point $i\in S$ is called a max point if 

i) $\exists k\in S$ such that $\mathcal{L}(i)\geq \mathcal {U}(k)$ and $\mathcal{L}(i)>\mathcal{L}(j)$, $\forall j\in (k,i)$, and

ii) $\exists k\in S$ such that $\mathcal{L}(i)\geq \mathcal {U}(k)$ and $\mathcal{L}(i)\geq\mathcal{L}(j)$, $\forall j\in (i,k)$.

A point $i\in S$ is called a min point if 

i) $\exists k\in S$ such that $\mathcal{U}(i)\leq \mathcal {L}(k)$ and $\mathcal{U}(i)<\mathcal{U}(j)$, $\forall j\in (k,i)$, and

ii) $\exists k\in S$ such that $\mathcal{U}(i)\leq \mathcal {L}(k)$ and $\mathcal{U}(i)\leq\mathcal{U}(j)$, $\forall j\in (i,k)$.
\end{defn}

\begin{lem} Between any two max (resp. min) points of $S$ there is a min (resp. max) point.\end{lem}

\textbf{Proof.} 
Let $m$ and $m'$ be two max points in $S$. Assume that $\mathcal{L}(m)\geq \mathcal{L}(m')$. There exists $k\in (m,m')$ such that $\mathcal{L}(m)\geq \mathcal{L}(m')\geq \mathcal{U}(k)$. Choose $k$ to be the first value where $\mathcal{U}$ attains a minimum on $(m,m')$ (the order here is the order induced by the cyclic order). Then $k$ is a min point. 
\QED

Let $m_1<m_2<\dots < m_l$ be the max and min points of $S$ and let $m_{l+1}=m_1$.

\begin{lem} Let $i\in\{1,2,\dots, l\}$. If $m_i$ is a max point then $\mathcal{U}(m_{i+1})\leq \mathcal{L}(m_i)$. If $m_i$ is a min point then $\mathcal{U}(m_{i})\leq \mathcal{L}(m_{i+1})$.\end{lem}
\textbf{Proof.} Let $m_i$ be a max point. Assume $\mathcal{U}(m_{i+1})> \mathcal{L}(m_i)$. Since $m_{i+1}$ is a min point, $\exists k\in(m_i,m_{i+1})$ such that $\mathcal{U}(m_{i+1})\leq \mathcal{L}(k)$. So $\mathcal{L}(m_i)<\mathcal{L}(k)<\mathcal{U}(k)$. Since $m_i$ is a max point, $\exists k'\in (m_i,k)$ such that $\mathcal{U}(k')\leq \mathcal{L}(m_i)<\mathcal{L}(k)$. Choose $k'$ to be the first point where $\mathcal{U}$ attains a minimum on $(m_i,k)$. Then $k'$ is a min point between $m_i$ and $m_{i+1}$, which is impossible. The proof of the second statement is similar.\QED

To simplify notation, fix $i\in\{1, 2, \dots, l\}$ and let $j=m_{i}$ and $k=m_{i+1}$. If $j$ is a max point then $L(t_k)<2\pi\mathcal{U}(k)\leq2\pi\mathcal{L}(j)<L(t_j)$, so the interval $(L(t_k), L(t_j))$ contains $\mathcal{L}(j)-\mathcal{U}(k)+1$ multiples of $2\pi$. Continuity of $L$ implies that each of these multiples of $2\pi$ is attained by $L$ at some value $t\in (t_j,t_k)$. Similarly, if $j$ is a min point, then each of the $\mathcal{L}(k)-\mathcal{U}(j)+1$ multiples of $2\pi$ in the interval $(L(t_j), L(t_k))$ is attained by $L$ at some value $t\in (t_j,t_k)$. Summing up over $i$, we get a value $V(f_0)$ which is a lower bound for the number of times $L$ crosses a multiple of $2\pi$ over the interval $[0, 2d-2]$. At each point $t$ where this happens, the corresponding function $f_t$ satisfies $f_t(r)=f_t(s)$. 

\begin{alg}\label{algorithm} The algorithm described above to compute $V(f_0)$ depends on the net of $f_0$ and not on $f_0$ itself, so we may think of it as accepting a net $\gamma$ and label its output $V(\gamma)$.
\end{alg}

\begin{figure}[p]
	\includegraphics[width=10cm]{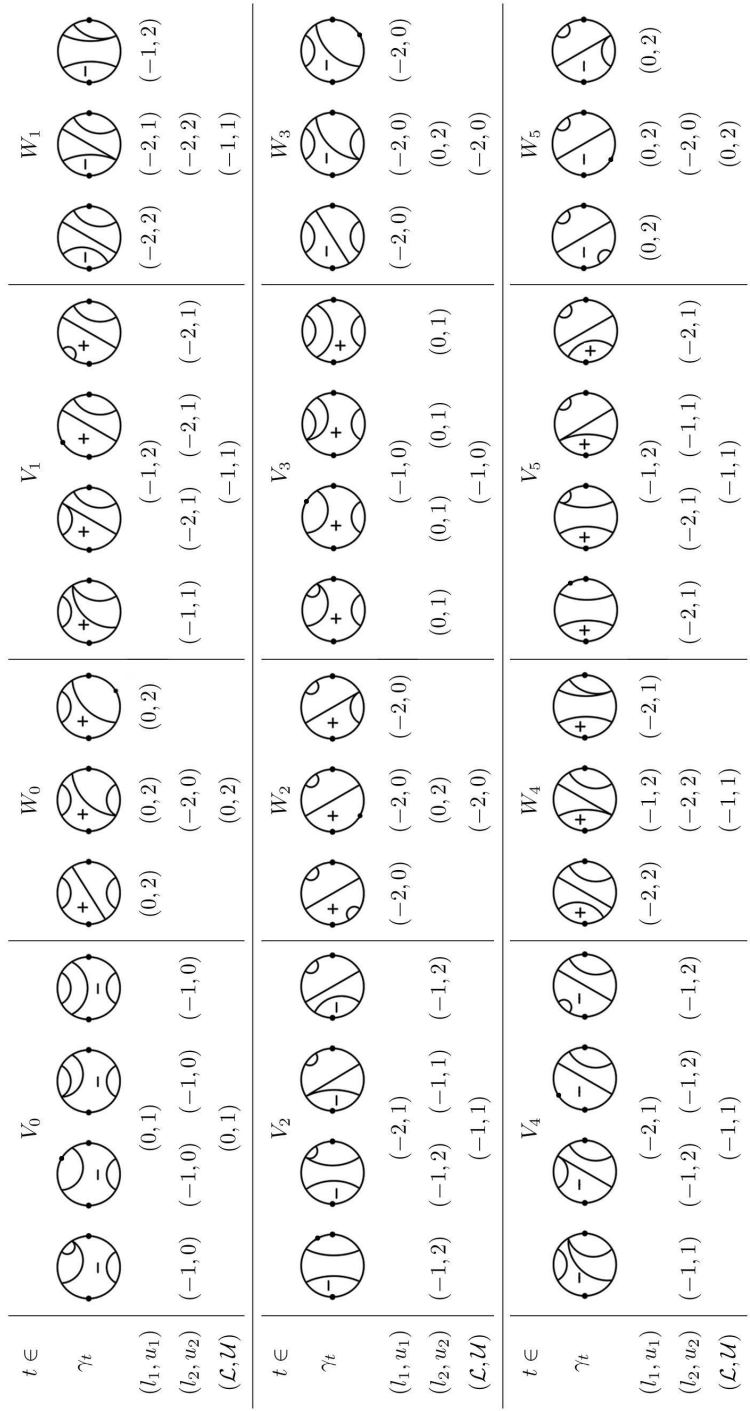}
	\caption{Some steps of Algorithm~\ref{algorithm} applied to a net with $d=4$ and $k=2$.}
	\label{minvar}
\end{figure}

\begin{ex}
Figure \ref{minvar} shows some steps of Algorithm~\ref{algorithm} applied to a net $\gamma$ with 8 vertices and $c=0$. 
The only min point of $\gamma$ corresponds to $W_2$ where $\mathcal{U}=0$, and the only max point corresponds to $W_5$ where $\mathcal{L}=0$. So $V(\gamma)=2$.
\end{ex}

Let $f_0^1,f_0^2,\ldots,f_0^{u_d}$ be nonequivalent functions in $\mathcal{S}_d$ having critical points at $1$  and at $y_1,y_2,\ldots, y_{2d-3}$. For each $i\in\{1,2,\ldots, u_d\}$ there exists a family $f_t^i$, $t\in[0,2d-2]$ satisfying 
properties (i)-(iii) of Lemma \ref{lifts}. 
By Theorem \ref{uniquec} and Notation \ref{rmk}, the map $t\longrightarrow [f_t^i]$ is unique. 
Let $\gamma_t^i=(\Gamma_t^i,y_1)$ be the net corresponding to $f_t^i$ and let $T_i$ be the period of the map $t\longrightarrow\gamma_t^i$. Uniqueness of the map $t\longrightarrow [f_t^i]$ implies that the map $t\longrightarrow\gamma_t^i$ is unique and that $T_i$ is well-defined. Let $C_i=\{\gamma_j^i:j=0,1,\ldots, T_i-1\}$. The sets $C_i$ are well-defined and partition the set $C$ of nets with $2d-2$ vertices. This follows from the fact that the map $t\longrightarrow\gamma_t^i$ is unique and its period is a factor of $2d-2$. Let $k$ be the number of distinct sets $C_i$. For $i\in\{1,2,\ldots,k\}$, let $N_{C_i}$ be the number of equivalence classes of functions $f$ in $\mathcal{S}_d$ whose net belongs to $C_i$ and which satisfy $f(r)=f(s)$. For $\gamma\in C$, let $V(\gamma)$ be the output of Algorithm~\ref{algorithm} when the input net is $\gamma$. 

\begin{thm} The number $N$ of classes $[f]\in\mathcal{S}_d$ having critical points at the points $y_1,y_2,\ldots, y_{2d-3}$ and satisfying $f(r)=f(s)$ is bounded below by $\sum_{\gamma\in C}\frac{\displaystyle V(\gamma)}{\displaystyle 2d-2} $.\end{thm} 

\textbf{Proof.} Since $C_i$, $i=1,2,\ldots,k$ partition $C$, $N=\sum_{i=1}^kN_{C_i}$. For each $i\in\{1,2,\ldots k\}$,
$$ N_{C_i}=\frac{\displaystyle T_i}{\displaystyle 2d-2}V({\gamma_0^i}).$$
Since $t\longrightarrow\gamma_t$ is periodic, $V(\gamma_j^i)=N_{C_i}$ for all $j\in\{0,1,\ldots,T_i-1\}$, hence,
$$ T_iN_{C_i}=\sum_{j=0}^{T_i-1} \frac{\displaystyle T_i}{\displaystyle 2d-2}V(\gamma_j^i).$$
Dividing both sides by $T_i$ and taking the sum over all classes $C_i$ we get
$$\sum_{i=1}^kN_{C_i}=\sum_{i=1}^k\sum_{j=0}^{T_i-1}\frac{\displaystyle V(\gamma_j^i)}{\displaystyle 2d-2}=\sum_{\gamma\in C}\frac{\displaystyle V(\gamma)}{\displaystyle 2d-2}. $$
\QED

The lower bounds computed by the algorithm for $4\leq d\leq 14$ and the corresponding values of $k$ appear in Table~\ref{lbtable}. The first row and column give the values of $d$ and $k$ respectively.

\begin{table}
	\centering
		\begin{tabular}[h]{|c|c|c|c|c|c|c|c|c|c|c|c|}
		\hline
		&4&5&6&7&8&9&10&11&12&13&14\\
		\hline
		1&1&4&14&48&165&572&2002&7072&25194&90440&326876\\
		\hline
		2&1&2&6&18&57&186&622&2120&7338&25724&91144\\
		\hline
		3&&4&12&36&113&366&1216&4122&14202&49592&175124\\
		\hline
		4&\multicolumn{2}{c|}\hfill&12&34&107&  {348}&  {1156}&3920&13514&47212&166788\\
		\hline
		5&\multicolumn{3}{c|}\hfill&36&  {115}&  {372}&  {1232}&4166&14326&49950&176178\\
		\hline
		6&\multicolumn{4}{c|}\hfill& 117& {370}&  {1232}&4164&14326&49920&175978\\
		\hline
		7&\multicolumn{5}{c|}\hfill&  {370}&  {1224}&4104&14024&48610&170606\\
		\hline
		8&\multicolumn{6}{c|}\hfill&   {1244}&4134&14176&49188&172660\\
		\hline
		9&\multicolumn{7}{c|}\hfill&4098&13948&48030&167690\\
		\hline
		10&\multicolumn{8}{c|}\hfill&14106&48348&169326\\
		\hline
		11&\multicolumn{9}{c|}\hfill&47904&166630\\
		\hline
		12&\multicolumn{10}{c|}\hfill&  168000\\
		\hline
		
		\end{tabular}
		\smallskip
		\caption{Lower bounds for $4\leq d\leq 14$.}
		\label{lbtable}
\end{table}
\section{Combinatorial interpretation for $k=1,2$}
In this section, we give a combinatorial interpretation of the algorithm in the cases $(r,s)$ contains exactly one or exactly two points of $D$. The first case agrees with the result obtained in section 4. 
For $t\in W_n$, $(l_1(t),u_1(t))\subseteq (c-u_2(t),c-l_2(t))$. For $t\in V_n$, $(l_1(t),u_1(t))\subseteq (c-u_2(t),c-l_2(t))$ except when $\theta(t)$ coincides with a fixed vertex $y$ and the two edges with endpoint $y$ have their other endpoint outside $(s,r)$. When $k=1$ and $t\in V_n$, there is only one fixed vertex outside $(s,r)$. So for $k=1$, it is enough to consider $(l_1,u_1)$.

\bigskip

\noindent\textbf{Case 1:} $(r,s)$ contains one point\\
Suppose there is only one fixed vertex $y$ in $(r,s)$. The Young tableaux  corresponding to a net with one fixed vertex in $(r,s)$ can have one or two squares. The only one having one square gives rise to a unique $(\mathcal{L},\mathcal{U})$-interval, namely $(-1,1)$. Each of the two tableaux with two squares gives rise to two intervals. These intervals are $(-1,0)$, $(0,1)$, $(-2,0)$ and $(0,2)$. All of these 5 intervals appear in Figure 6. 

The three sequences $(0,1) (-1,1) (-1,0)$, $ (0,1) (-1,1) (-2,0)$, $ (0,2) (-1,1) (-1,0)$, and their reversals cannot occur in the sequence of $(\mathcal{L},\mathcal{U})$-intervals of any net with $d>2$. So $(-1,0)$ and $(0,1)$ do not yield any max or min points. The $(\mathcal{L},\mathcal{U})$-interval $(0,2)$ gives a max point if and only if it is preceded by $(-2,0) (-1,1)$. In this case the net corresponding to $(-1,1)$ has the property that the moving vertex is outside $(r,s)$ and the fixed vertex in $(r,s)$ is not connected to any of its neighboring vertices. The same holds when $(-2,0)$ is a min point. Also, if the net has the property that the moving vertex is outside $(r,s)$ and the fixed vertex in $(r,s)$ is not connected to any of its neighboring vertices then its $(\mathcal{L},\mathcal{U})$-interval must be preceded by $(0,2)$ and followed by $(-2,0)$ or vice versa. 

So given a net $\Gamma$ with $2d-2$ vertices, $V(\Gamma)$ is precisely the number of vertices $y$ of $\Gamma$ not connected to any of their neighboring vertices. The lower bound for a given $d$ and $k=1$ is the number of distinct nets $(\Gamma,y)$ such that $\Gamma$ has $2d-2$ vertices and $y$ is not connected to any of its neighboring vertices. This agrees with the result obtained in section 4.

\begin{figure}[h]\begin{center}
	\includegraphics{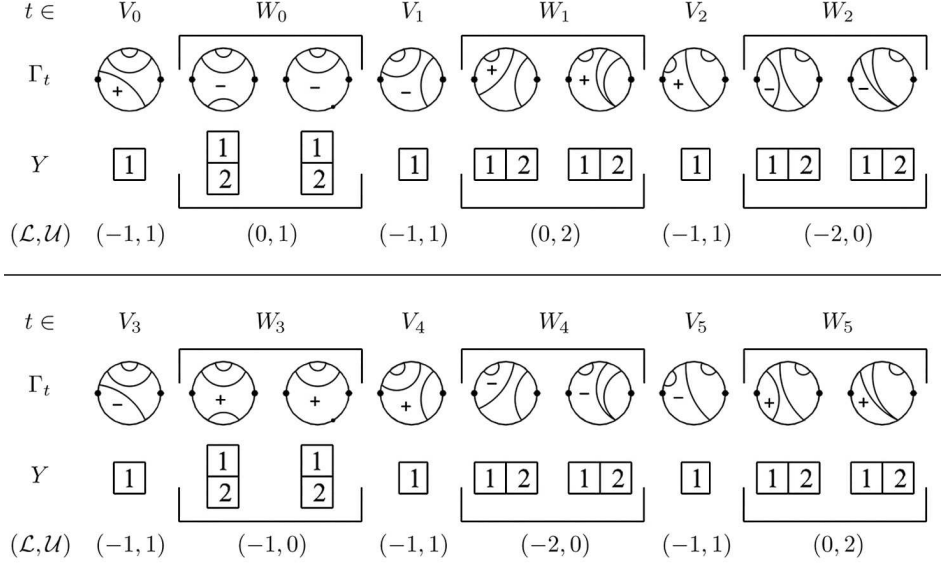}
 	\caption{Algorithm~\ref{algorithm} applied to a net with $d=4$ and $k=1$.}
	\label{minvar1}\end{center}
\end{figure}

\begin{ex}
Figure \ref{minvar1} shows some steps of Algorithm~\ref{algorithm} applied to a net $\gamma$ with 6 vertices and $k=1$. The degenerate nets with degenerate vertex outside $(r,s)$ are not shown. Below each net is its Young tableau corresponding to $(r,s)$. The last row gives the $(\mathcal{L},\mathcal{U})$-interval for each $V_i$ and $W_i$, $i=0,1,\dots, 5$. The only min point of $\gamma$ corresponds to $W_2$ where $\mathcal{U}=0$, and the only max point corresponds to $W_5$ where $\mathcal{L}=0$. So $V(\gamma)=2$, which is the number of vertices of $\Gamma$ not connected to any of the neighboring vertices.
\end{ex}
\bigskip
\noindent\textbf{Case 2:} $(r,s)$ contains two points

Suppose there are two fixed vertices $y_1$ and $y_2$ in $(r,s)$. All possible $(\mathcal{L},\mathcal{U})$-intervals appear in Figure 5. Any max point must correspond to an $(\mathcal{L},\mathcal{U})$-interval with $\mathcal{L}=0$. More specifically, any max point must correspond to an $(\mathcal{L},\mathcal{U})$-interval equal to $(0,2)$ since the interval $(0,1)$ must always be preceded by $(0,2)$. Similarly, a min point  must correspond to an $(\mathcal{L},\mathcal{U})$-interval equal to $(-2,0)$. 

Any interval giving rise to a max (resp. min) point must be immediately followed by the interval $(0,1)$ (resp. $(-1,0)$). The interval $(0,1)$ (resp. $(-1,0)$) is obtained when the moving critical point is outside $(r,s)$, the two critical points in $(r,s)$ are connected by a chord and the function is orientation reversing (resp. preserving) in a neighborhood of $r$. 

So the min/max points correspond to chords in the net connecting adjacent critical points and separated by an odd number of critical points. The lower bound is $\mathcal{N}_{2d-2}$, the number  of nets $\gamma'=(\Gamma',y_1)$ with $2d-2$ vertices satisfying $y_1y_2$ is a chord of $\Gamma'$ and the next chord connecting two consecutive vertices is of the form $y_{2k}y_{2k+1}$ with $1<k<d-1$. 

Counting the number of such nets for a fixed $k$ is equivalent to finding the number $N_{k}$ of nets $\gamma=(\Gamma,v_1)$ having $2d-6$ vertices $v_1,v_2,\ldots,v_{2d-6}$ such that if $v_iv_{i+1}$ is a chord of $\Gamma$ then $i\geq 2k-3$. Let $S$ be the set of all nondegenerate nets $\gamma$ having $2d-6$ vertices. Given a net $\gamma$, let $C_\gamma$ be the set of chords of $\gamma$. 

$N_2=u_{d-2}$ since, when $k=2$ there is no condition on the chords of $\gamma$.

$N_3$  is the number of nets $\gamma$ such that $v_1v_2$ and $v_2v_3$ are not chords in $\gamma$. So $N_3=u_{d-2}-2u_{d-3}$. 

$N_4$ is the number of nets $\gamma$ such that $v_iv_{i+1}$ is not a chord in $\gamma$ for $i=1,2,3,4$. 
\begin{eqnarray}
N_4 &=& |S|
     -\sum\limits_{i=1}^4 | \{\gamma\in S : v_iv_{i+1}\in C_{\gamma}\}|\cr
    && + |\big\{\gamma\in S : \{v_1v_2,v_3v_4\} \mbox{ or } \{v_1v_2,v_4v_5\} \mbox{ or } \{v_2v_3,v_4v_5\} \subset C_{\gamma}\big\}|\cr\nonumber
    &=& u_{d-2}-4u_{d-3}+3u_{d-4}.\cr\nonumber
\end{eqnarray}
Similar computations yield $N_{5}=u_{d-2}-6u_{d-3}+10u_{d-4}-4u_{d-5}$ and $N_{6}=u_{d-2}-8u_{d-3}+21u_{d-4}-20u_{d-5}+5u_{d-6}$, etc. 
Using these values, we can compute $\mathcal{N}_{2d-2}$. 
The values of $\mathcal{N}_{2d-2}$ for some small values of $d$ are
\begin{eqnarray}
\mathcal{N}_6&=&N_2=u_2=1\cr\nonumber
\mathcal{N}_8&=&N_2+N_3=2u_3-2u_2=2\cr\nonumber
\mathcal{N}_{10}&=&N_2+N_3+N_4=3u_4-6u_3+3u_2=6\cr\nonumber
\mathcal{N}_{12}&=&N_2+N_3+N_4+N_{5}=4u_5-12u_4+13u_3-4u_2=18\cr\nonumber
\mathcal{N}_{14}&=&N_2+N_3+N_4+N_{5}+N_{6}=5u_6-20u_5+34u_4-24u_3+5u_2=57.\cr\nonumber
\end{eqnarray}

\bibliography{bibliography}
\bibliographystyle{plain}
\end{document}